\documentclass{article}
\usepackage{arxiv}

\usepackage[utf8]{inputenc} 
\usepackage[T1]{fontenc}    
\usepackage{hyperref}       
\usepackage{url}            
\usepackage{booktabs}       
\usepackage{amsfonts}       
\usepackage{nicefrac}       
\usepackage{microtype}  

\usepackage{amsthm}
\usepackage{amsmath}
\usepackage{graphicx}
\usepackage{amssymb}
\usepackage{color}
\usepackage{tikz}
\usepackage{hyperref}
\usepackage{algorithm}
\usepackage{algorithmic}

\usepackage{multirow}
\usetikzlibrary{arrows,shapes,chains}

	\title{An adaptive homotopy method for computing bifurcations of nonlinear parametric systems}

\author{
	Wenrui Hao \\
	Department of Mathematics \\
	The Pennsylvania State University, University Park, PA 16802, USA.\\
	\texttt{wxh64@psu.edu} \\
	\And
	Chunyue Zheng \\
	Department of Mathematics\\
	Pennsylvania State University, University Park, PA 16802\\
	The Pennsylvania State University, University Park, PA 16802, USA.\\
	\texttt{cmz5199@psu.edu} \\
}
\begin{document}
\maketitle

\begin{abstract}
	In this paper, we present an adaptive step-size homotopy tracking method for computing bifurcation points of nonlinear systems. There are four components in this new method: 1) an adaptive tracking technique is developed near bifurcation points; 2) an inflation technique is backed up when the adaptive tracking fails; 3) Puiseux series interpolation is used to compute bifurcation points; and 4)
the tangent cone structure of the bifurcation point is approximated numerically to compute solutions on different branches.
Various
numerical examples of nonlinear systems are given to illustrate the efficiency of this new approach. This new adaptive homotopy tracking method is also applied to a system of nonlinear PDEs and shows robustness and efficiency for large-scale nonlinear discretized systems.
\end{abstract}


\section{Introduction}
\label{sec:into}

Many mathematical models of natural phenomena, e.g.,
biology \cite{HHHS}, physics
\cite{HNS,HNS1} and materials science \cite{HouLowengrub}, involve
systems of nonlinear equations
\cite{FH,HCF,HF,HHHS}. From a mathematical point of view, studies
of these nonlinear equations can be formulated numerically and
theoretically to focus on solution structures such as
bifurcations \cite{Rhein4,Rhein1}. Theories and numerical
methods have contributed to a better understanding of these
solution structures, in which the bifurcation between solutions and
parameters is the central question \cite{haber1,Strogatz}. Although
theory helps us to understand the solution structures in many
cases \cite{FHB,FHB1}, the in-depth and more quantitative
study of these problems often requires large-scale simulations
to numerically compute bifurcations.
A bifurcation occurs in a  nonlinear parametric system when the parameter change causes the solution structure to change. There are many types of bifurcations, such as saddle-node bifurcation, transcritical bifurcation, pitchfork bifurcation, and Hopf bifurcation with different theoretical classifications \cite{khalil2002nonlinear}. However computing these different bifurcation points numerically brings the same numerical challenge. In specific, this corresponds to the real part of an eigenvalue of the Jacobian passing through zero and causes numerical challenges for Newton's and Newton-like methods \cite{dayton2005computing,zeng2005computing,zeng2004algorithm}.
Therefore, efficient numerical methods for
computing bifurcations of large-scale systems are keys to
understanding these systems.

The homotopy continuation method
\cite{LiSauer,LiZeng,MorganSommese1,MorganSommese} has been
successfully used to compute bifurcations and structural stabilities for studying parametric problems. Recently, several numerical
methods based on homotopy continuation methods have been developed for computing
bifurcation points of
nonlinear PDEs \cite{HHHLSZ,HHHS}. These numerical methods have also been
applied to hyperbolic conservation laws \cite{HHSSXZ}, physical
systems \cite{HNS,HNS1} and some more complex free boundary problems
arising from biology \cite{HCF,HF}. However, the computational cost becomes more expensive and the efficiency becomes low when they are
applied to large-scale systems.
Therefore,
an efficient homotopy continuation method for computing
bifurcation is
needed to deeply study the large-scale nonlinear systems. In this paper, we will present an efficient adaptive homotopy tracking method that integrates numerical methods from numerical algebraic geometry and scientific computing so that we can apply this efficient method to compute bifurcation points of large-scale nonlinear systems such as discretized systems arising from nonlinear PDEs.

\section{Homotopy Continuation Method}
\label{sec:main}
In this section, we will first give an overview of the homotopy continuation method. 	Generally speaking, a  nonlinear parametric system is written as $\mathbf{F}:
\mathbb{R}^n\times\mathbb{R}\rightarrow\mathbb{R}^n,$
\begin{equation}\label{Sys}
\mathbf{F}(\mathbf{u},p)=\mathbf{0},
\end{equation}
where $p$ is a parameter and $\mathbf{u}$ is the variable vector \cite{BHS,MSW} that depends on the parameter $p$, namely, $\mathbf{u}=\mathbf{u}(p)$.
We want to start with solutions that are easy to find (e.g.,
radially symmetric solutions in nonlinear PDEs \cite{HHHLSZ}) in order to compute
the bifurcation points where the other more interesting solutions
come from  (e.g., non-radial solution \cite{HHHLSZ}).
For this parametric system, the standard homotopy continuation method \cite{CHB,WBM} uses a
predictor/corrector method to track the solution $\mathbf{u}$ as
the parameter $p$ varies. Basic prediction and correction are both
accomplished by considering a local model via its Taylor expansion:
$$
\mathbf{F}(\mathbf{u}+\Delta \mathbf{u},p+\Delta
p)=\mathbf{F}(\mathbf{u},p)+\mathbf{F}_\mathbf{u}(\mathbf{u},p)\Delta
\mathbf{u}+\mathbf{F}_p(\mathbf{u},p)\Delta p+\hbox{Higher-Order Terms},
$$
where $\mathbf{F}_\mathbf{u}=\partial \mathbf{F}/\partial \mathbf{u}$ is the $n\times
n$ Jacobian matrix and $\mathbf{F}_p=\partial \mathbf{F} /\partial p$ has size
$n\times1$.

\subsection{Predictor-Corrector Method}
\begin{figure}[ht]
	\centering
	\includegraphics[width=0.5\linewidth]{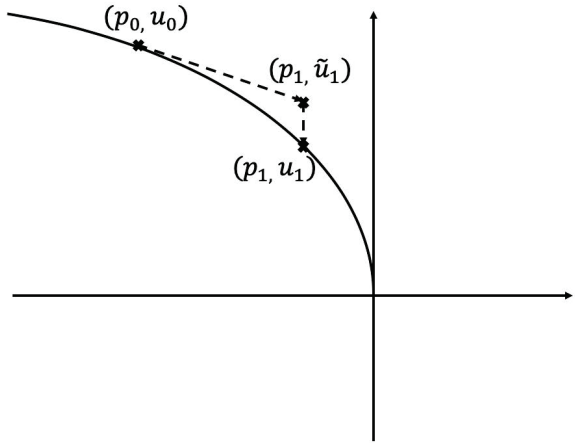}
	\caption{An illustration of the predictor-corrector Method.}
	\label{Fig:PC}
\end{figure}

The Predictor-Corrector method consists of two parts: the first one is the predictor step which gives a prediction of $\Delta \mathbf{u}$ for any given $\Delta p$ based on numerical methods for solving ordinary differential equation such as  Euler method, the secant predictor method, and etc (see \cite{allgower2003introduction} for more details); the second one is the corrector method which refines the predicted solution based on numerical methods for solving nonlinear systems such as Newton's method, conjugate gradient methods and etc (see \cite{allgower2003introduction} for more details). In this section, we will use the Euler predictor and the Newton corrector to illustrate the idea of the predictor-corrector procedure. Other predictor-corrector method can be found in \cite{allgower2003introduction}.
{Given a solution $(\mathbf{u}_0,p_0)$ on the path, that
	is, $\mathbf{F}(\mathbf{u}_0,p_0)=0$, we plan to compute a solution at
	$p_1=p_0+\Delta p$ by setting $\mathbf{F}(\mathbf{u}_0+\Delta \mathbf{u},p_0+\Delta
	p)=0$. First we make an Euler predictor step, solving the first-order terms $
	\mathbf{F}_\mathbf{u}(\mathbf{u}_0,p_0)\Delta \mathbf{u}=-\mathbf{F}_p(\mathbf{u}_0,p_0)\Delta
	p$, and then letting $\tilde{\mathbf{u}}_1=\mathbf{u}_0+\Delta\mathbf{u}$. On the other hand, when $\|\mathbf{F}(\tilde{\mathbf{u}}_1,p_1)\|$ is not sufficiently small, one may fix $p_1$ to be constant by setting
	$\Delta p=0$ and solving the following equation by using the Newton corrector: $
	\mathbf{F}_\mathbf{u}(\tilde{\mathbf{u}}_1,p_1)\Delta
	\mathbf{u}=-\mathbf{F}(\tilde{\mathbf{u}}_1,p_1).$ Repeat this corrector step until $\|\mathbf{F}(\tilde{\mathbf{u}}_1,p_1)\|$ is smaller than the chosen tolerance criterion, then we can get $\mathbf{u}_1=\tilde{\mathbf{u}}_1+\Delta\mathbf{u}$ and $(\mathbf{u}_1,p_1)$ is on the path (see an illustration in Fig. \ref{Fig:PC}).
}

\subsection{The Step-Size Control}
The main concern of a numerical path-tracking algorithm is deciding
which of these methods
to use next and how large of a step-size $\Delta p$ to use in
the predictor \cite{bates2006bertini,deuflhard2011newton}.
A trial-and-error approach for the step-size control is
used for homotopy continuation tracking: shorten the step-size upon failure and lengthen it upon repeated successes \cite{BHSWbook,SW}.
This trial-and-error approach can be computationally expensive and can lack efficiency when systems are not well-conditioned,  since the step-size becomes very small.
Moreover, in the path tracking process, at some critical
points, the ill-conditioned Jacobian matrix $F_\mathbf{u}$  often
causes trouble either in the prediction or in the correction process.
Various computational techniques, such as pseudo-arclength continuation, Gauss-Newton continuation, and other adaptive step-size strategies \cite{deuflhard2011newton}, have been developed to handle this difficulty. For instance, the path tracking may encounter no difficulty at a turning point if the pseudo-arclength continuation is adopted. However, bifurcations of large-scale nonlinear systems are usually complex (more than turning points) and need a more sophisticated numerical method to compute.

\section{Adaptive Homotopy Tracking with Bifurcation Detection (AHTBD)}
To overcome this difficulty, an adaptive homotopy tracker is proposed to reduce the computational cost. The basic idea of this adaptive homotopy tracker is to solve the step-size simultaneously when we track the nonlinear system.
For any given step-size $h$, we start with  a point on the solution path, denoted by $(\mathbf{u}_0,p_0)$, and want to find the next point to satisfy the following augmented system:
\begin{eqnarray}\tilde{\mathbf{F}}(\mathbf{u},p)=\left(
\begin{array}{c}
\mathbf{F}(\mathbf{u},p) \\ g \mathbf{v}^T(\mathbf{u}-\mathbf{u}_0)(1-s)+s(p-p_0)-h
\end{array}
\right),\label{ASC}\end{eqnarray} where
$g=\mathrm{sign}(-\mathbf{v}^T\mathbf{F}_\mathbf{u}(\mathbf{u}_0,p_0)^{-1}\mathbf{F}_p(\mathbf{u}_0,p_0))/\| \mathbf{F}_\mathbf{u}(\tilde{\mathbf{u}},\tilde{p})^{-1}\mathbf{F}_p(\tilde{\mathbf{u}},\tilde{p})\|$, $s=\Big|\frac{\lambda_{min}}{\tilde{\lambda}_{min}}\Big|$,  $\lambda_{min}$ is the real part of the minimum eigenvalue of
$\mathbf{F}_\mathbf{u}$ at $(\mathbf{u}_0,p_0)$, and $\mathbf{v}$ is the corresponding eigenvector. Here $(\tilde{\mathbf{u}},\tilde{p})$ is a generic point (i.e., randomly choosing $\tilde{p}$ to compute $\tilde{\mathbf{u}}$) \cite{BHSWbook,SW} and $\tilde{\lambda}_{min}$ is the real part of the minimum eigenvalue of
$\mathbf{F}_\mathbf{u}$ at $\tilde{p}$.
Thus the next point on the path $(\mathbf{u},p)$ is computed by solving the new augmented system $\tilde{F}$ with an adaptive step-size. In particular,
when  the tracking parameter $p$
is close to a bifurcation point, $\lambda_{min}$ is very small, and $s$ approaches zero, we then have
$g\mathbf{v}^T(\mathbf{u}-\mathbf{u}_0)=h$ instead of $p-p_0=h$ which
means that we change the tracking parameter from $p$ to
$\mathbf{v}^T\mathbf{u}$; when $p_0$ is a generic point, namely, the original system is well-conditioned, we have $s$  be close to $1$ and then $p=p_0+h$ which is the ``initial" target for the next point.
Moreover,  this adaptive homotopy tracking process, whose pseudocode is outlined in {\bf Algorithm \ref{alg1}}, employs the Newton-Krylov method  to solve the augmented nonlinear system.
%
%

\begin{algorithm}[H]
	\caption{The pseudocode of the adaptive tracking algorithm.}\label{alg1}
	\begin{algorithmic}
		\STATE \textbf{Input: }{A step-size $h$,  a start point $(\mathbf{u}_0,p_0)$, and an ending parameter $p_e$}.
		\STATE \textbf{Output: }{A solution sequence on the path $(\mathbf{u}_i,p_i)_{i=1}^N$}. 
		\STATE Set $i=0$;
		\WHILE {$(p-p_0)(p-p_e)\leq 0$}
		\STATE	Compute the minimum eigenvalue of $\mathbf{F}_\mathbf{u}(\mathbf{u}_i,p_i)$ and the corresponding eigenvector, $\mathbf{v}$;
		\STATE	Solve the augmented system \eqref{ASC} and denote the solution as $(\mathbf{u}_{i+1},p_{i+1})$;
		\STATE	Set $i=i+1$;
		\ENDWHILE
	\end{algorithmic}
\end{algorithm}

\noindent	{\bf Remark 1:} The augmented system (\ref{ASC}) does not bring new singularities. In other words, if the original system is full rank, then the augmented system must be full rank. In fact, if $\mathbf{F}_\mathbf{u}$ is not singular, the Jacobian matrix of the augmented system (\ref{ASC}) can be written as
$$\left(
\begin{array}{cc}
\mathbf{F}_\mathbf{u} & \mathbf{F}_p\\
g\mathbf{v}^T(1-s) &s
\end{array}
\right) = \left(
\begin{array}{cc}
I & 0 \\
g\mathbf{v}^T(1-s)\mathbf{F}_\mathbf{u}^{-1} & I
\end{array}
\right)\left(
\begin{array}{cc}
\mathbf{F}_\mathbf{u} & \mathbf{F}_p \\
0 & s-g\mathbf{v}^T(1-s)\mathbf{F}_\mathbf{u}^{-1}\mathbf{F}_p
\end{array}
\right). $$
If the original system has full rank, namely, $s\neq 0$,
then we have $s-(1-s)g\mathbf{v}^T\mathbf{F}^{-1}_\mathbf{u}\mathbf{F}_p\ne0$, which implies that the augmented system (\ref{ASC}) also has full rank. On the other hand, if $F_u$ is singular, the Jacobian matrix of the augmented system could be non-singular.

\noindent	{\bf Remark 2:} The parameter tracking direction is the same as $h$. In fact, by solving
\begin{equation*}
\left(
\begin{array}{cc}
\mathbf{F}_\mathbf{u} & \mathbf{F}_p\\
g\mathbf{v}^T(1-s) &s
\end{array}
\right)\left(\begin{array}{c}\Delta\mathbf{u}\\\Delta p\end{array}\right)=\left(\begin{array}{c}0\\h\end{array}\right),
\end{equation*}
we have
$$
\Delta p=\frac{h}{s-(1-s)g\mathbf{v}^T\mathbf{F}^{-1}_\mathbf{u}\mathbf{F}_p}.
$$
Noticing the definition of $g$, we have $s-(1-s)g\mathbf{v}^T\mathbf{F}^{-1}_\mathbf{u}\mathbf{F}_p>0$ if $s\ne0$, which implies that $\Delta p$ has the same sign as $h$.	\subsection{Inflation Process}
When the Jacobian matrix of the augmented system is ill-conditioned, the adaptive path tracking algorithm
based on Newton's method is no longer satisfactory since it may converge slowly or even diverge. Once such a circumstance occurs, the deflation technique has been proposed to overcome this difficulty \cite{leykin2006newton,hauenstein2013isosingular}. However, the deflated system is double the size of the original nonlinear system, and sometimes
even higher order derivatives
need to be taken into consideration \cite{leykin2006newton}. Therefore this technique is hard to apply for large-scale systems. In order to track large-scale systems, we need a different strategy, an inflation process. The motivation of the inflation technique is based on iterative methods for the ill-conditioned symmetric positive definite matrices. Let us consider a simple example with $(A+\epsilon I)x=b$ ($A$ and $b$ are shown below), and apply the Gauss-Seidel method with stopping criteria $\|Ax^k-b\|\leq 10^{-8}$ and $x^0=b$.  Eq. (\ref{Tab1}) shows the number of iterations for different value of $\epsilon$: the number of iterations increases dramatically from $18$ to $54,470$ when the matrix is ill-conditioned; the number of iterations drops to 2 when the matrix is singular. Therefore iterative methods usually are effective for a singular system, but time-consuming for a nearly singular system (see \cite{lee2007robust} for more theoretical results).
\begin{equation}\label{Tab1}
\centering
A=\begin{bmatrix}
1      & -1 & 0 \\
-1     & 2 & -1 \\
0 & -1 & 1
\end{bmatrix}, b=\begin{bmatrix}
-1\\
-1 \\
2
\end{bmatrix}\in R(A).\quad\quad\begin{tabular}{|c|c|}
\hline
$\epsilon$ & \# of of iteration \\
\hline
1 & 18\\
$10^{-1}$&100\\
$10^{-2}$ & 852 \\
$10^{-3}$ & 6,982\\
$10^{-4}$ & 54,470\\
$0$ & 2\\\hline
\end{tabular}
\end{equation}
Based on this motivation, we will inflate the nearly singular system to a singular system. More specifically,
for a bifurcation point $p^*$, the system $\mathbf{F}(\mathbf{u}^*,p^*)$ is
singular. By denoting $J$ the Jacobian $\mathbf{F}_\mathbf{u}(\mathbf{u},p)$,
we know that $J$ is ill-conditioned if $p$ is close to $p^*$ so that
Newton's method becomes difficult to converge. By decomposing
$\Delta \mathbf{u}$ as $\Delta \mathbf{u}=\widetilde{\Delta \mathbf{u}}+\alpha
\mathbf{v}$, then we solve the following inflated system instead of
$\mathbf{F}_\mathbf{u}(\mathbf{u},p)\Delta
\mathbf{u}=-\mathbf{F}(\mathbf{u},p)$:
\begin{eqnarray}\left(
\begin{array}{cc}
J^TJ & J^TJ\mathbf{v}\\
\mathbf{v}^T J^TJ&\lambda_{min}
\end{array}
\right) \left(
\begin{array}{c}
\widetilde{\Delta \mathbf{u}} \\ \alpha
\end{array}
\right)=-\left(
\begin{array}{c}
J^TF(\mathbf{u},p) \\ \mathbf{v}^TJ^TF(\mathbf{u},p)
\end{array}
\right).  \label{Aug}\end{eqnarray}
Here $\lambda_{min}$ is the eigenvalue of $J^TJ$ with the minimum norm and $\mathbf{v}$ is the corresponding eigenvector. We use $J^TJ$ instead of $J$ to make sure the coefficient matrix is symmetric positive semi-definite in order to guarantee the convergence of this inflation technique  \cite{lee2007robust}. In fact, for any $a\in\mathbb{R}^{n\times1}, b\in\mathbb{R}$, we have
\begin{equation}
\begin{aligned}
(a^T, b)\left(
\begin{array}{cc}
J^TJ & J^TJ\mathbf{v}\\
\mathbf{v}^T J^TJ&\lambda_{min}
\end{array}
\right)\left(\begin{array}{c}a\\b\end{array}\right)&=a^TJ^TJa+b\mathbf{v}^T J^TJa+a^TJ^TJ\mathbf{v}b+\lambda_{min}b^2\\
&=a^TJ^TJa+2\lambda_{min}ba^T\mathbf{v}+\lambda_{min}b^2\\
&\ge\lambda_{min}|a|^2-2\lambda_{min}|b||a||\mathbf{v}|+\lambda_{min}b^2\\
&\ge\lambda_{min}(|a|-|b|)^2,
\end{aligned}
\end{equation}
which implies that the matrix in (\ref{Aug}) is symmetric positive semi-definite.
Therefore linear iterative solvers
such as Gauss-Seidel or GMRES \cite{XuSiam,XuLong} converge very quickly for solving the singular inflated  system   (\ref{Aug}) \cite{lee2007robust}.


\noindent {\bf Remark:} 	Since $(\mathbf{v}^T,-1)^T$ is in the kernel of \eqref{Aug}, we have a family of solutions $(\widetilde{\Delta \mathbf{u}}+k\mathbf{v},\alpha-k)$   for \eqref{Aug}, $\forall k$, for any given solution pair
$(\widetilde{\Delta \mathbf{u}},\alpha)$. However $\Delta \mathbf{u}$ is unique for any $k$ by the definition.

\subsection{Puiseux Series Extrapolation}
The power series endgame has been successfully used to
handle the singularity in NAG \cite{BHS,MSW} for polynomial systems. This endgame technique
is only used for homotopy tracking very near  $t=0$, but cannot
handle the bifurcation point during the tracking. In this paper,
we will develop a new numerical method based on the Puiseux Series Expansion (PSE) to approximate the
bifurcation point and the solution at the
bifurcation point when the nonlinear system is polynomial. The idea is to use the eigenvalue of the Jacobian matrix to interpolate the solution near the bifurcation point. In particular,
at the bifurcation point, the Jacobian
$\mathbf{F}_\mathbf{u}$ has an eigenvalue with zero real part, say
$p_b$, and several branches can come together at
$(\mathbf{u}_b,p_b)$. We denote $\lambda=\min_i|real(\lambda_i)|$,
where $\lambda_i$ is the eigenvalue of
$\mathbf{F}_\mathbf{u}(\mathbf{u},p)$ for any given
$(\mathbf{u},p)$. Then according to the classical Puiseux's theorem (Chapter 7 in \cite{fischer2001plane} \& Corollary A.3.3 in \cite{SW})
we use a Puiseux series expansion to
approximate $(\mathbf{u},p)$ in a neighborhood of
$(\mathbf{u}_b,p_b)$, called the PSE
operating zone. Thus the following formulation is given by
\begin{eqnarray}
\mathbf{u}(\lambda)=\mathbf{u}_b+\sum_{j=1}^\infty\mathbf{a}_i\lambda^{j/c_1}\hbox{~and~}p(\lambda)=p_b+\sum_{j=1}^\infty
b_i\lambda^{j/c_2}, \label{PSE}\end{eqnarray} where $c_1$ and $c_2$ are
the winding numbers for path $\mathbf{u}(\lambda)$ and
$p(\lambda)$, respectively.
Computing the winding numbers
$c_1$ and $c_2$ requires more advanced computational techniques in
NAG \cite{BHSWbook,huber1998polyhedral,SW} but can not be applied directly for large-scale nonlinear systems, e.g., the discretized
polynomial systems of nonlinear PDEs.
Thus in our algorithm, we make several guesses at
$c_1$ and $c_2$ to get the close connection to the curvature of the
paths.

Moreover, we also need to compute  leading terms of the PSE, namely, $w=\min\{j|\mathbf{a}_j\neq 0\}$ and  $q=\min\{j|b_j\neq 0\}$. Then (\ref{PSE}) is rewritten as
\begin{eqnarray} \mathbf{u}(\lambda)=\mathbf{u}_b+\lambda^{w/c_1}(\mathbf{a}_w+\sum_{j=w+1}^\infty\mathbf{a}_i\lambda^{j/c_1})\hbox{~and~}p(\lambda)=p_b+\lambda^{q/c_2}(b_q+\sum_{j=q+1}^\infty	b_i\lambda^{j/c_2}). \end{eqnarray}     We will show the procedure how to estimate $q/c_2$, which can be extended to estimate $w/c_1$ as well: for any constant $k_1$ and $k_2$, we have \begin{equation*}
\begin{aligned}
p(k_1\lambda)&=p_b+k_1^{q/c_2}\lambda^{q/c_2}(b_q+\sum_{j=q+1}^\infty	b_i(k_1\lambda)^{j/c_2}),\\
p(k_2\lambda)&=p_b+k_2^{q/c_2}\lambda^{q/c_2}(b_q+\sum_{j=q+1}^\infty	b_i(k_2\lambda)^{j/c_2}).
\end{aligned}
\end{equation*}
When $\lambda$ is small and $k_1<1,k_2<1$, we have
$$\frac{1-k_1^{q/c_2}}{1-k_2^{q/c_2}}\approx\frac{p(\lambda)-p(k_1\lambda)}{p(\lambda)-p(k_2\lambda)}.$$
Thus an approximation of $q/c_2$ is obtained by solving the following nonlinear equation:
$$f(x):= 1-k_1^x -m(1-k_2^x)=0,$$
where $m=\frac{p(\lambda)-p(k_1\lambda)}{p(\lambda)-p(k_2\lambda)}$.
{For estimating $w/c_1$,  we multiply a random vector, $\mathbf{\alpha}$, namely, using $\alpha^T\mathbf{u}(k_1\lambda)$ and $\alpha^T\mathbf{u}(k_2\lambda)$ to repeat the above procedure.}	In summary,	the algorithm for computing the bifurcation point based on the PSE is as follows:


\begin{algorithm}[H]
	\caption{Implementing PSE }
	\begin{algorithmic}
		\STATE Given a sequence of points on the branch, say $(\mathbf{u}^n,p^n,\lambda^n)_{n=1}^N$.
		\WHILE {$|\lambda|<Tol$}
		\STATE	Estimate the value of $w/c_1$ and $q/c_2$ by solving the nonlinear equation $f(x)=0$;
		\FOR { $c_i = 1:M$}
		\STATE Use the first $N-1$ points to approximate the Puiseux series;
		\STATE Apply these approximations to extrapolate $(\mathbf{u}^N,p^N)$ at $\lambda^N$;
		\ENDFOR
		\STATE	Determine the best value of $c_i$  by  choosing the nearest extrapolating point on the paths at $\lambda=\lambda^N$;
		\STATE	Use the Puiseux series to approximate $(\mathbf{u}_b,p_b)$ at $\lambda=0$;
		\IF {$\|(\mathbf{u}_b,p_b)\|<Tol$}
		\STATE Break;
		\ELSE 
		\STATE Set $\lambda=\frac{\lambda_N}{2}$, generate a new point $(\mathbf{u}^{N+1},p^{N+1})$, and update the sequence of points;
		\ENDIF
		\ENDWHILE
	\end{algorithmic}
\end{algorithm}

\begin{figure}[th]
	\centering
	\includegraphics[width=0.45\linewidth]{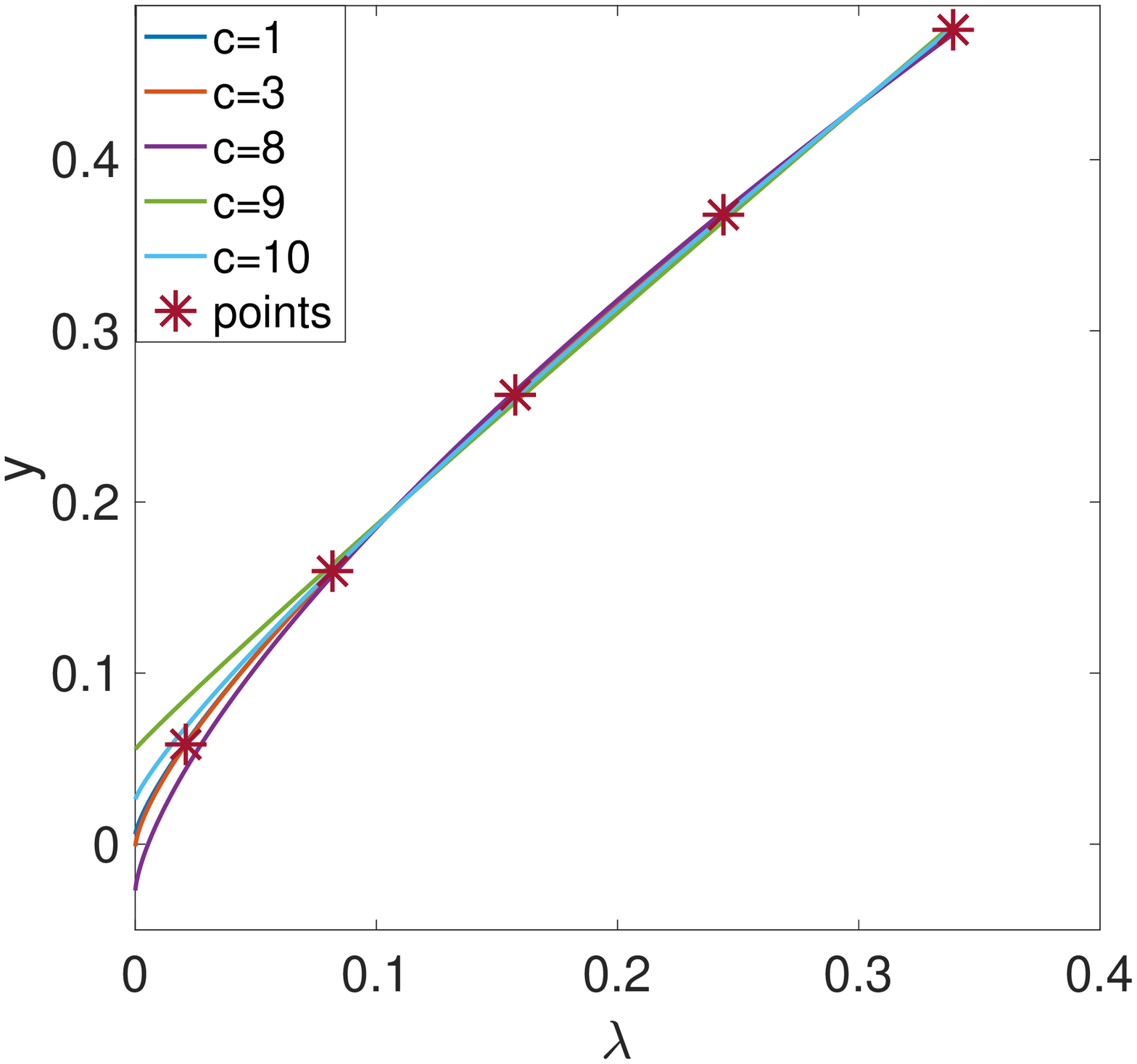}
	\includegraphics[width=0.45\linewidth]{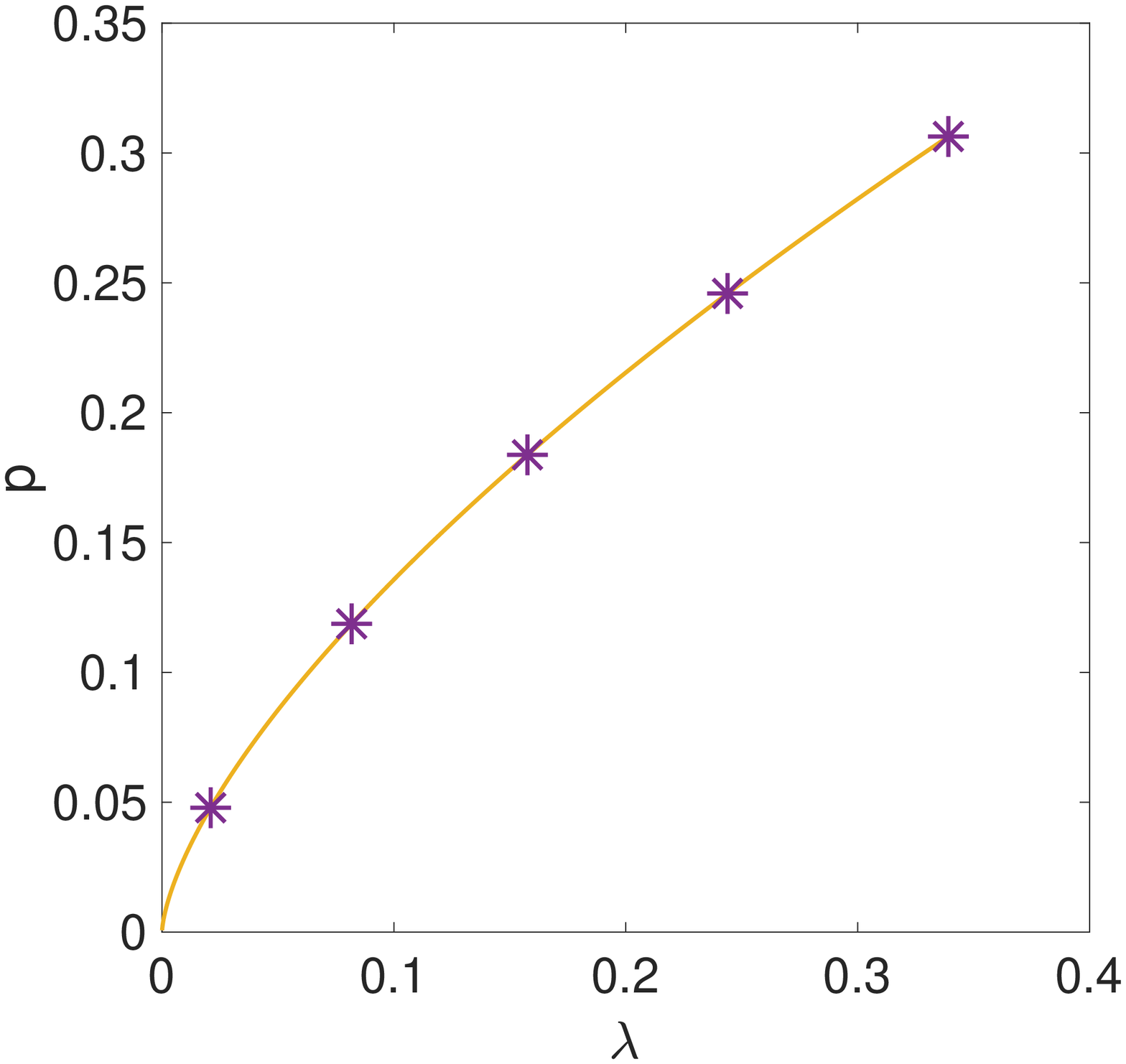}
	\caption{The PSE interpolation in the illustrated example (\ref{ILEX}). The left part shows solution trajectories of $y$ with respect to $\lambda$ for  different $c_1$; the right part shows parameter $p$ with respect to $\lambda$. }
	\label{fig:PSE}
\end{figure}

%
\noindent	{\bf An illustrated example: } We will use the following example to illustrate this PSE interpolation process:
\begin{equation}
F(\mathbf{u},p)=\left(\begin{matrix}
x^2-p^2\\
(x+y)^2-p^3
\end{matrix}\right).\label{ILEX}
\end{equation}
In this example, exact solutions of one branch are
$$
x=-\Big(\frac{1}{2}\Big)^{2/3}\lambda^{2/3},\ y=\Big(\frac{1}{2}\Big)^{1/3}\lambda^{2/3}+\frac{1}{2}\lambda \hbox{~and~} \ p=\Big(\frac{1}{2}\Big)^{2/3}\lambda^{2/3},
$$
where $\lambda$ is the minimum eigenvalue of the Jacobian matrix. By taking $\lambda=2$, we have our initial point $x_0=-1$, $y_0=2$, and $p_0=1$. By taking $h=-0.1$, we collect five points on this solution path  shown in Fig. \ref{fig:PSE}. Four of them are used to compute coefficients of the Puiseux series, the other one is to determine the winding numbers $c_1$ and $c_2$. Fig. \ref{fig:PSE} shows different solution trajectories by using PSE interpolation for different $c_1$. Then $c_1=3$ is the best approximation for $x$, $y$. In fact, since $p$ is a monomial of $\lambda$, when using a different winding number $c_2$, the ratio $q/c_2$ is the same. Then the approximated bifurcation point becomes $x= -3.2\times10^{-5},~ y= 1.1\times10^{-4}$, and $p=3.2\times10^{-5}.$

\subsection{Tangent Cone} After computing the bifurcation point,  the tangent cone of the bifurcation point needs to be computed in order to track along different branches { by using the Lyapunov-Schmidt reduction \cite{buffoni2003analytic,chicone1994lyapunov,HHHLSZ}}.
The tangent cone $T_*$ and the Jacobian matrix $J_*$ at the bifurcation point have the following relationship
$$
T_*\subseteq null(J_*),
$$
which implies that the tangent cone is contained in the tangent space at a bifurcation although the tangent cone and tangent space are equal at a generic  point.
Then the null space of the Jacobian is computed to obtain the tangent cone at a bifurcation  by using the Taylor expansion of the nonlinear system $\mathbf{F}$ in the null space of $J_*$. We will illustrate the procedure of computing the tangent cone by assuming that the dimension of the null space of the $J_*$ is $n-1$. Let's denote the Jacobian $J_{\mathbf{u}}$ and the derivative $J_p$ with respect to $p $ at $(\mathbf{u}_0,p_0)$ as $A:=[J_\mathbf{u},J_p]\in R^{n\times(n+1)}$. Then we have \[\begin{bmatrix}
\mathbf{Q}_1      & \mathbf{Q}_2\\
q_1     & q_2
\end{bmatrix} =null(A), \hbox{~where~} \mathbf{Q}_i\in R^{n\times 1}  \hbox{~and~} q_i \hbox{~is a scalar}.\] Similarly, $\Lambda\in R^{1\times n}=null(A^T)$. Thus we assume that
\[\Delta \mathbf{u}=a_1 \mathbf{Q}_1+a_2 \mathbf{Q}_2 \hbox{~and~} \Delta p=a_1 q_1+a_2 q_2,\] where $a_i$ needs to be determined. We construct the following single polynomial $g(a_1,a_2)$
\[g(\mathbf{a})=\Lambda^TF(\mathbf{u}_0+a_1\mathbf{Q}_1+a_2\mathbf{Q}_2,p_0+a_1q_1+a_2q_2).\]
By using Taylor expansion at $(0,0)$, we have
\[g(\mathbf{a})\approx g(0,0)+\mathbf{a}^T\frac{\partial g}{\partial \mathbf{a}}(0,0)+\mathbf{a}^TH(0,0)\mathbf{a},\]
where $H(0,0)$ is  the Hessian matrix of $g$ at $(0,0)$. Then $\mathbf{a}$ stratifies the following system:
\begin{equation*}
\begin{aligned}
\mathbf{a}^TH(\mathbf{0})\mathbf{a}&=0\\
{a}_1q_1+a_2q_2&=\Delta_p.
\end{aligned}
\end{equation*}
If the tangent cone has a more complex structure (such as when the dimension of the null space of the Jacobian is more than 1), we need to introduce more variables $a_i$ and more derivatives to determine the tangent cone.

\begin{figure}[th]
	\centering
	\tikzstyle{process} = [rectangle, minimum width = 2cm, minimum height = 1cm, text centered, draw = black]
	\tikzstyle{decision} = [diamond, aspect = 3, text centered, draw=black]
	\begin{tikzpicture}[font=\sffamily]
	\node (O) [process] {\small Given a starting point and a step-size};
	\node (D) [decision, below =1 cm of O] {\small Is the augmented system nearly singular?};
	\node (X1) [process,below of=D, yshift=-1cm,, right=4cm] {\small Inflation process};
	\node (X2) [process,below of=D, yshift=-1cm, left=3cm] {\small Adaptive homotopy tracking};
	\node (Y) [process,below =6cm  of O]  {\small Compute the bifurcation point by using the PSE extrapolation};
	\node (Z) [process,below =1 cm of Y] {\small Compute the tangent cone to obtain the local bifurcation structure};
	
	\draw [semithick,-] (O) --(D);
	\draw [semithick,->] (D) --node[below]{yes} (X1);
	\draw [semithick,->] (D) --node[below]{no} (X2);
	\draw [semithick,->] (X2) |-(D);
	\draw [semithick,->] (X1) |-(D);
	\draw [semithick,->] (X2) |-(Y);
	\draw [semithick,->] (X1) |-(Y);
	\draw [semithick,->] (Y) -- (Z);
	\draw [semithick,-] (Z) -| (7,-1);
	\draw [semithick,->] (7,-1) -- (0,-1);
	\end{tikzpicture}
	\caption{The flow chart of the AHTBD method.}
	\label{alg}
\end{figure}
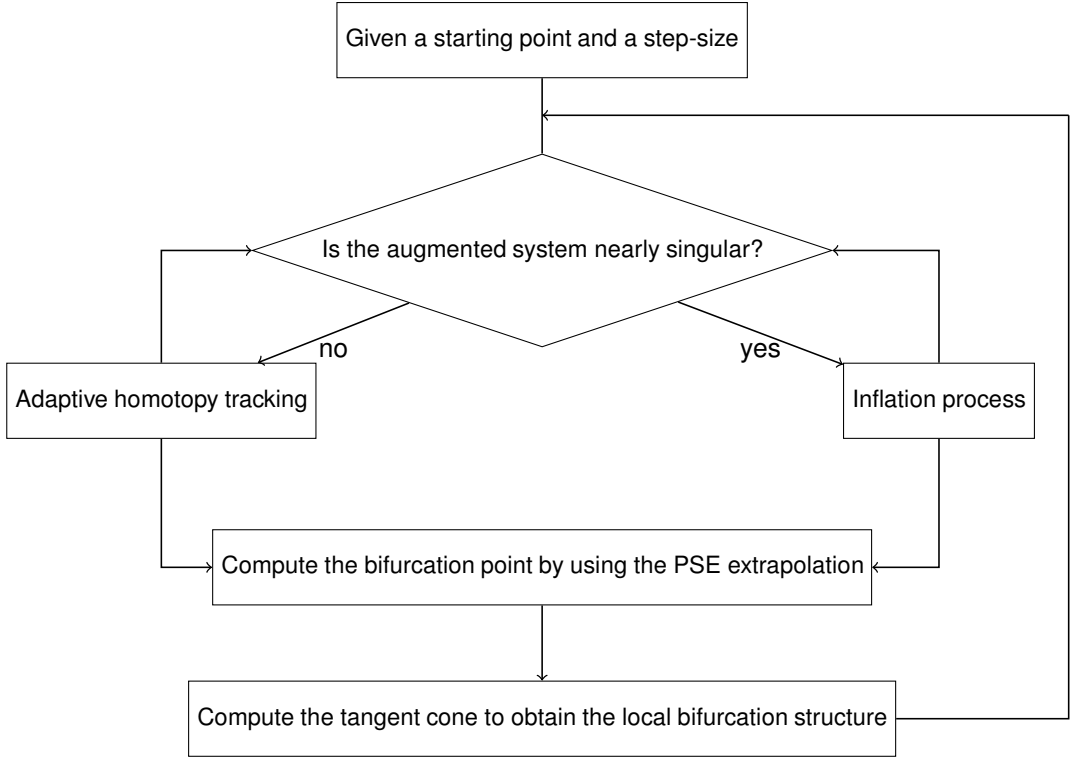

Therefore, we summarize the AHTBD method as follows and outline the flow chart in Fig. \ref{alg}:
\begin{enumerate}
	\item For a given initial point $(\mathbf{u},p)$ on a solution path and a maximum step-size, solve the augmented system (\ref{ASC}) to track along the path;
	
	\item If  the augmented system (\ref{ASC}) becomes {ill-conditioned}, the inflation process is introduced;
	
	\item Near the bifurcation point,  the PSE interpolation is used to approximate the bifurcation point;
	
	\item At the bifurcation point, the tangent cone is computed to determine the different tracking solution branches, and then repeat the first step for each path.
	
\end{enumerate}

\section{Numerical Results}
In this section, we apply the AHTBD method to several examples, ranging from a single equation to a system of nonlinear PDEs, to show its efficiency. Both the AHTBD method and the traditional homotopy tracking method are implemented and compared on Matlab.  The traditional homotopy tracking has been implemented in various packages such as Bertini \cite{bates2006bertini}, HOM4PS \cite{LLT}, PHCpack \cite{verschelde1999algorithm} and others to handle the bifurcations. {Among these existing software, Bertini has more freedom to compute the  bifurcations due to the adaptive multi-precision path tracking \cite{bates2008adaptive} and the parallel endgame \cite{bates2011parallel}.} To fairly compare the AHTBD method with the traditional homotopy tracking, we will implement both methods on Matlab.

\subsection{An example with a turning point}
\label{subsec: compare} 
\begin{figure}[th]
	\centering
	\includegraphics[width=0.45\linewidth]{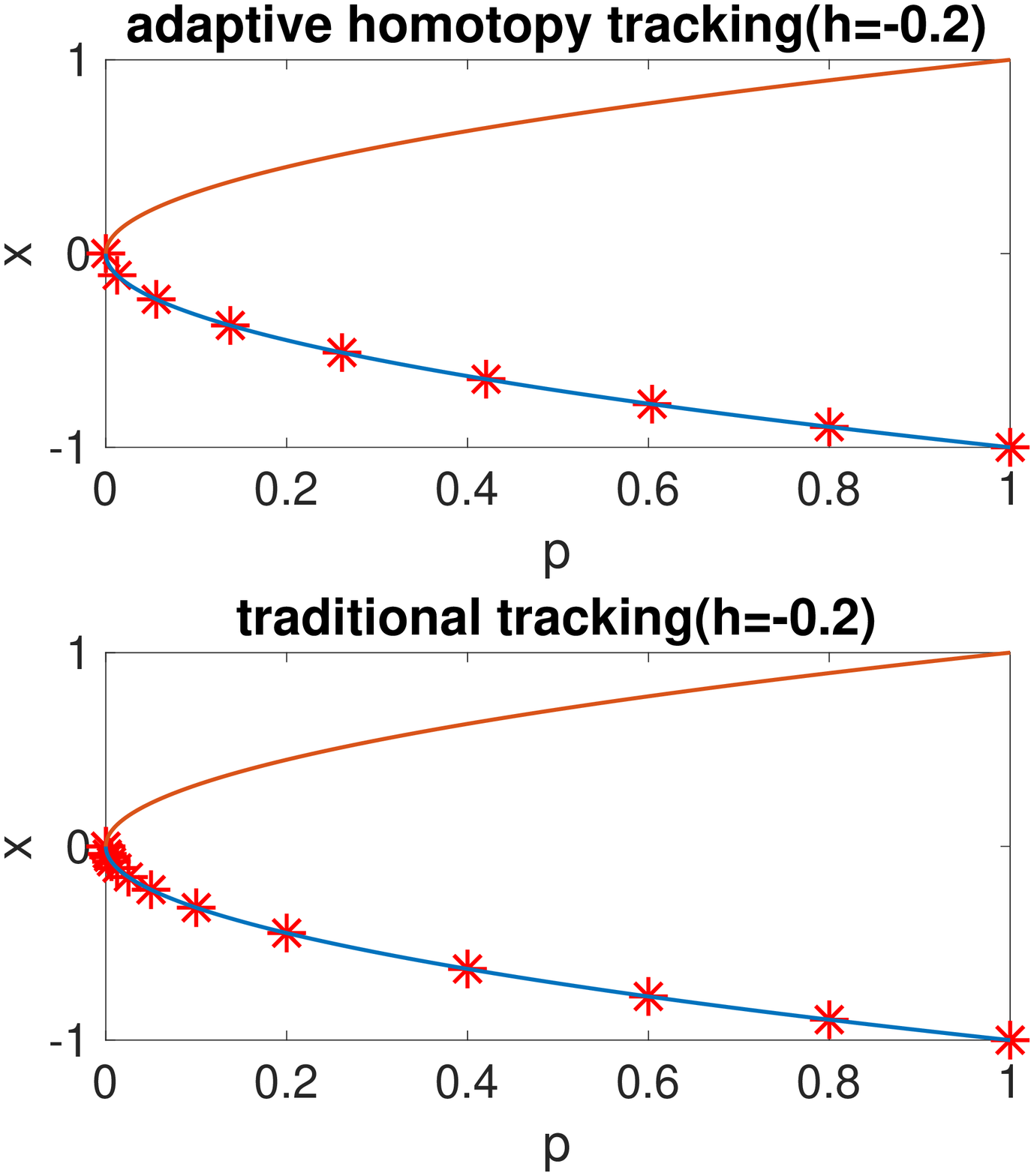}
	\includegraphics[width=0.45\linewidth]{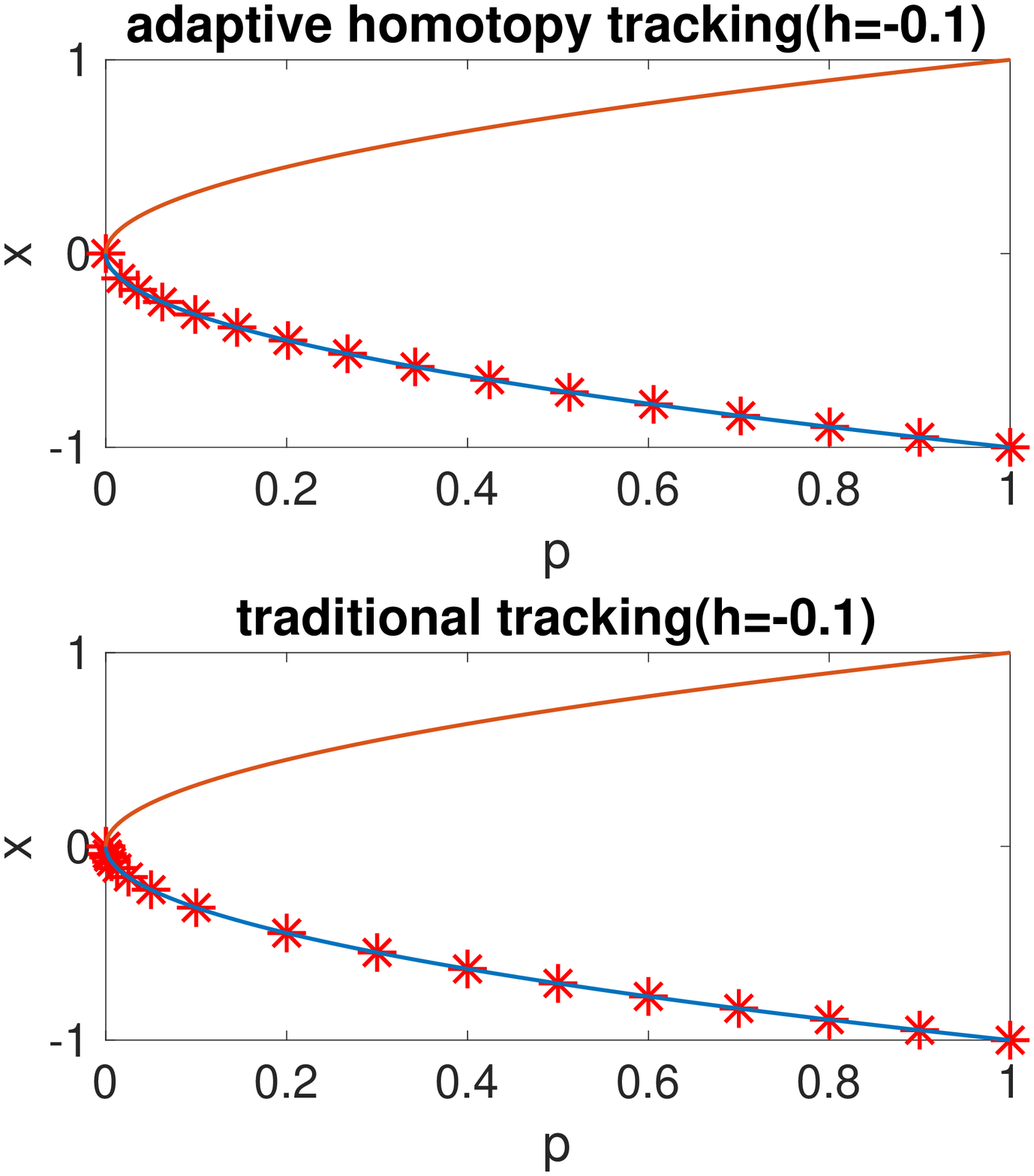}
	\caption{Comparisons between adaptive (upper) and traditional (lower) homotopy tracking methods. The plot of  $x$ v.s. $p$ is illustrated for $h=-0.2$ (left) and $h=-0.1$ (right).}
	\label{Fig:compare_h}
\end{figure}

Our first example is used to test the efficiency of adaptive homotopy tracker by considering the following system:
\begin{equation}
F(\mathbf{u},p)=\left(\begin{matrix}
x^2-p\\
x^2-2y^2+p
\end{matrix}\right),\label{ex11}
\end{equation}
where $\mathbf{u}=(x,y)^T$ is the variable while $p$ is the parameter. The analytical solution is $x^2=y^2=p$ which has a turning point when $p=0$. This example is used to illustrate the efficiency of the adaptive homotopy tracker for computing the bifurcation point. We choose  $\mathbf{u}_0=(-1,1)$ and $p_0=1$ as our initial tracking point and compare the adaptive homotopy tracker of the AHTBD method and the traditional homotopy tracker with different step-sizes ($h=-0.1$ and $h=-0.2$). Table \ref{ex1Table} and Fig.~\ref{Fig:compare_h} show that the adaptive homotopy tracker takes fewer steps to get to the bifurcation point. In particular, when the initial step-size $h$ becomes larger, the efficiency of the adaptive homotopy tracker is more obvious. The traditional homotopy method finds the bifurcation by halving the step-size  with less accuracy (around $10^{-4}$), while the adaptive homotopy tracker approximates the bifurcation point by doing the PSE extrapolation with higher accuracy  (around $10^{-6}$).

\begin{table}
	\begin{center}
		\begin{tabular}{|c|c|c|c|c|c|}
			\hline
			\multicolumn{2}{|c|}{\multirow{2}{*}{}}
			&\multicolumn{2}{|c|}{h=-0.1}&\multicolumn{2}{|c|}{h=-0.2}\\\cline{3-6}
			\multicolumn{2}{|c|}{\multirow{2}{*}{}}	& traditional & adaptive & traditional & adaptive\\
			\hline
			\multicolumn{2}{|c|}{\# of steps}  & 19 & 16  &15  & 9  \\ \hline
			{\multirow{3}{*}{bifurcation}}& x& $3\times10^{-4}$  & $-1.695\times10^{-5}$  &  $3\times10^{-4}$  & $-1.1\times 10^{-16}$\\\cline{2-6}
			&y& $-3\times10^{-4}$  & $1.695\times10^{-5}$  &  $-3\times10^{-4}$  & $1.1\times 10^{-16}$\\\cline{2-6}
			&p	& $-5.7\times10^{-7}$  & $5.5\times10^{-6}$&  $-5.7\times10^{-7}$  & $3.8\times10^{-11}$\\\hline
		\end{tabular}
		\caption{Comparisons between adaptive and traditional  homotopy tracking methods for (\ref{ex11}).}
		\label{ex1Table}
	\end{center}
\end{table}

\subsection{Examples with complex bifurcation structures}
In this subsection, we will use the AHTBD method to compute several examples with complex bifurcation structures; namely,  the bifurcation point is computed first by using the adaptive homotopy tracker, and then the tangent cone algorithm is used to obtain different solution branches.

\begin{figure}[th]
	\centering
	\includegraphics[width=0.5\linewidth]{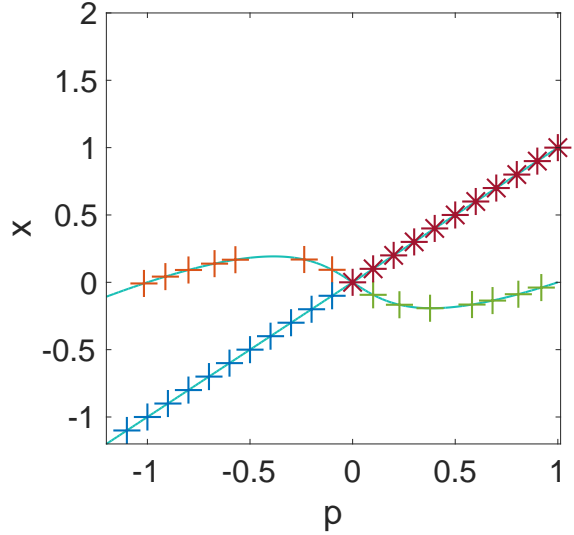}
	\caption{Local bifurcation diagram of (\ref{ex21}): starting from the lower branch (blue points), we compute the bifurcation point first by using the PSE interpolation and then compute the tangent cone to obtain the other solution branches (green, red, and orange points).}
	\label{Fig:ex21}
\end{figure}

\noindent{\bf Example 1:}  Given
\begin{equation}
F(x,p)=(x-p)^4+(x-p)(x+p),\label{ex21}
\end{equation}
we have a bifurcation point at $p=0$. In order to compute the local bifurcation diagram at $p=0$, we start from a point $x=1$ and $p=1$ to track along a solution path with the step-size $h=-0.1$. {When it is close to the bifurcation, namely, $\lambda_{min}<0.1$, we use the PSE to approximate the bifurcation point. Afterwards the tangent cone is computed: since the Jacobian $F_x$ and the derivative $F_p$ are both 0, the null space for $A=[F_x,F_p]$ is $\mathrm{span}\{(0,1)^T,(1,0)^T \}$ and the null space of $A^T$ is $\mathrm{span}\{1\}$. Then two tangent directions are obtained, $(1,1)^T$ and $(-1,1)^T$}. By setting different step-sizes, for example $h=\pm0.1$, and choosing a tangent direction,  we obtain a solution on each branch. Starting from this point, the  adaptive homotopy algorithm is employed to continue tracking (see Fig. \ref{Fig:ex21}).

\begin{figure}[th]
	\centering
	\includegraphics[width=0.45\linewidth]{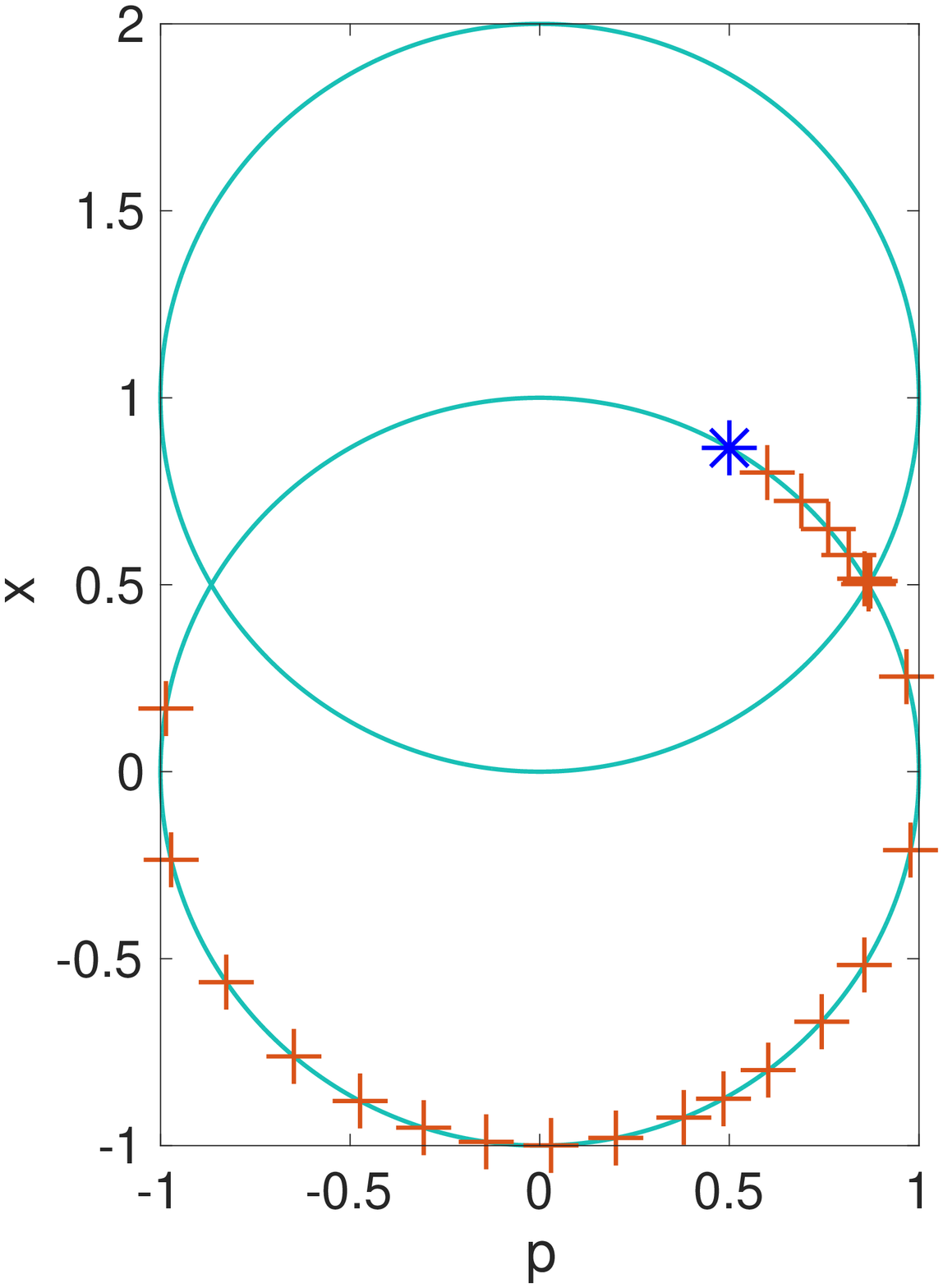}
	\includegraphics[width=0.45\linewidth]{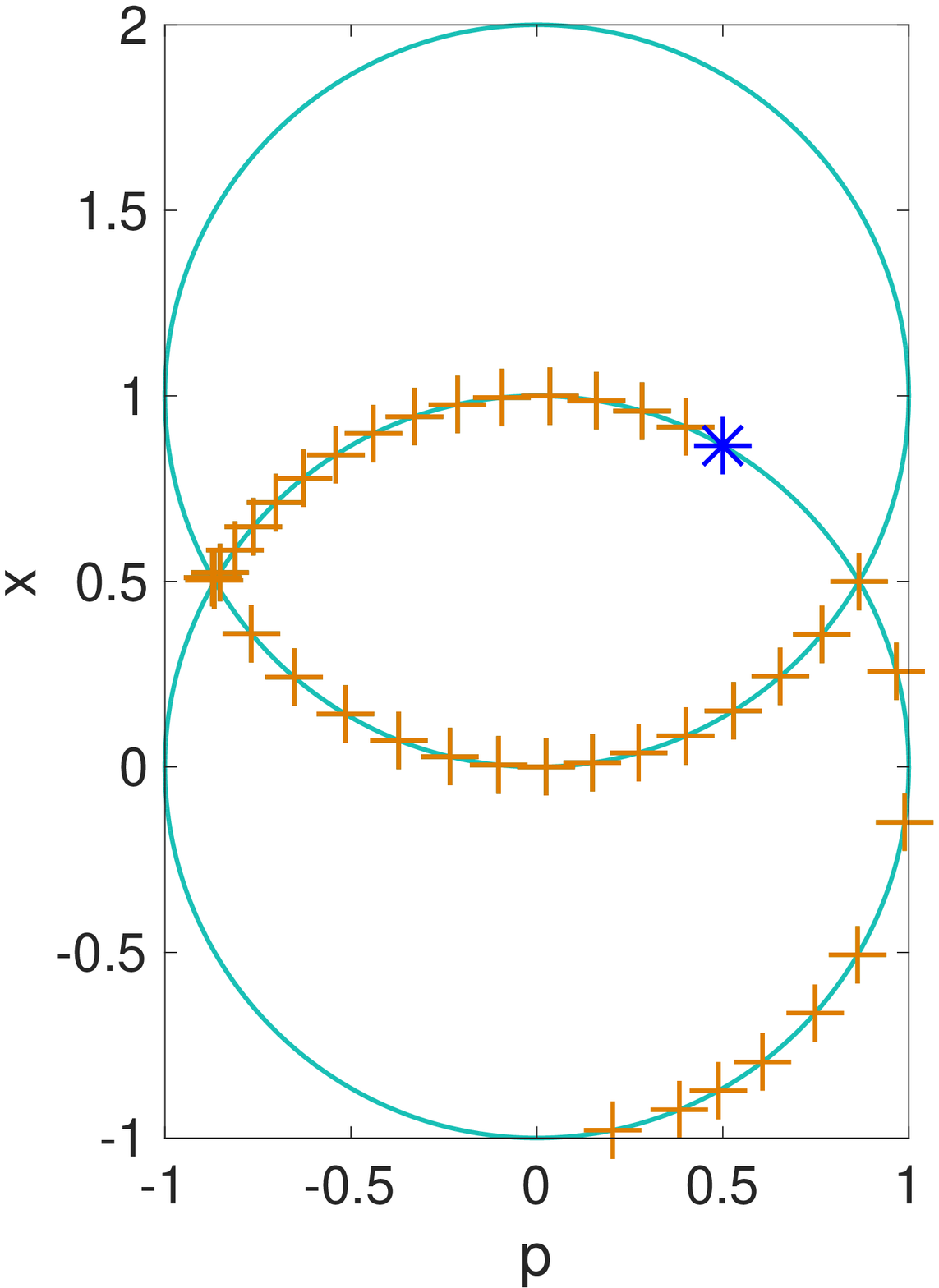}
	\caption{Local bifurcation diagram of (\ref{ex22}). The AHTBD method is used to track from the blue point to the left and right directions.}\label{Fig:ex22}
\end{figure}

\noindent	{\bf Example 2:} The following equation represents two intersecting circles that imply complex bifurcation structures shown in Fig. \ref{Fig:ex22}:
\begin{equation}
F(x,p)=(x^2+p^2-1)((x-1)^2+p^2-1).\label{ex22}
\end{equation}
We start to track along a solution path from point $(\frac{1}{2},\frac{\sqrt{3}}{2})$ with different tracking directions (blue point in Fig. \ref{Fig:ex22}). Fig. \ref{Fig:ex22} shows the AHTBD tracking process with the step-size $|h|=0.1$. It is clearly seen that the tracking is almost uniform even though there are two bifurcation points. Table \ref{tab:ex2_compare} shows the comparison between the AHTBD and traditional homotopy methods when the tracking starts at point $(\frac{1}{2},\frac{\sqrt{3}}{2})$ and ends when reaching or passing the turning point where $|p|=1$. The two tables have the same starting point, while the tracking direction is different. Although the traditional homotopy method may have higher accuracy for the bifurcation point, it takes many more steps to reach the end point than the AHTBD method. Moreover, the AHTBD method can pass the turning point easily (see Table \ref{tab:ex2_compare} for $h=-0.1$), while the traditional method stagnates at the turning point.

\begin{table}
	\begin{center}
		\begin{tabular}{|c|c|c|c|c|c|}
			\hline
			\multicolumn{2}{|c|}{\multirow{2}{*}{}}
			&\multicolumn{2}{|c|}{h=0.1}&\multicolumn{2}{|c|}{h=0.05}\\\cline{3-6}
			\multicolumn{2}{|c|}{\multirow{2}{*}{}}	& Trial-and-error   & AHTBD & Trial-and-error   & AHTBD\\
			\hline
			\multicolumn{2}{|c|}{\# of steps}  & 50 & 10  &60  & 17 \\ \hline
			{\multirow{2}{*}{bifurcation}}& x& $0.5002$  & $0.5020$  &  $0.5002$  & $0.5033$\\\cline{2-6}
			&p	& $0.8659$  & $0.8671$&  $0.8659$  & $0.8702$\\\hline
			{\multirow{2}{*}{endpoint}}& x& $0.0128$  & $-0.2097$  &  $0.0128$  & $-0.1020$\\\cline{2-6}
			&p	& $0.9999$  & $0.9778$&  $0.9999$  & $0.9948$\\\hline
		\end{tabular}
		\begin{tabular}{|c|c|c|c|c|c|}
			\hline
			\multicolumn{2}{|c|}{\multirow{2}{*}{}}
			&\multicolumn{2}{|c|}{h=-0.1}&\multicolumn{2}{|c|}{h=-0.05}\\\cline{3-6}
			\multicolumn{2}{|c|}{\multirow{2}{*}{}}	& Trial-and-error   & AHTBD & Trial-and-error   & AHTBD\\
			\hline
			\multicolumn{2}{|c|}{\# of steps}  & 62 & 19  &82  & 32  \\ \hline
			{\multirow{2}{*}{bifurcation}}& x& $0.5002$  & $0.5026$  &  $0.5002$  & $0.4831$\\\cline{2-6}
			&p	& $-0.8659$  & $-0.8675$&  $-0.8659$  & $-0.8756$\\\hline
			{\multirow{2}{*}{endpoint}}& x& $0.0080$   & $-0.1637$  &  $0.0080$  & $0.0054$\\\cline{2-6}
			&p	& $-1.0000$ & $-0.9865$&  $1.0000$  & $-1.0000$\\\hline
		\end{tabular}
		\caption{Comparisons between AHTBD and traditional trial-and-error  tracking  methods along the branches shown in Fig. \ref{Fig:ex22} with different step-sizes for $h$.}
		\label{tab:ex2_compare}
	\end{center}
\end{table}

\begin{figure}[th]
	\centering
	\includegraphics[width=0.45\linewidth]{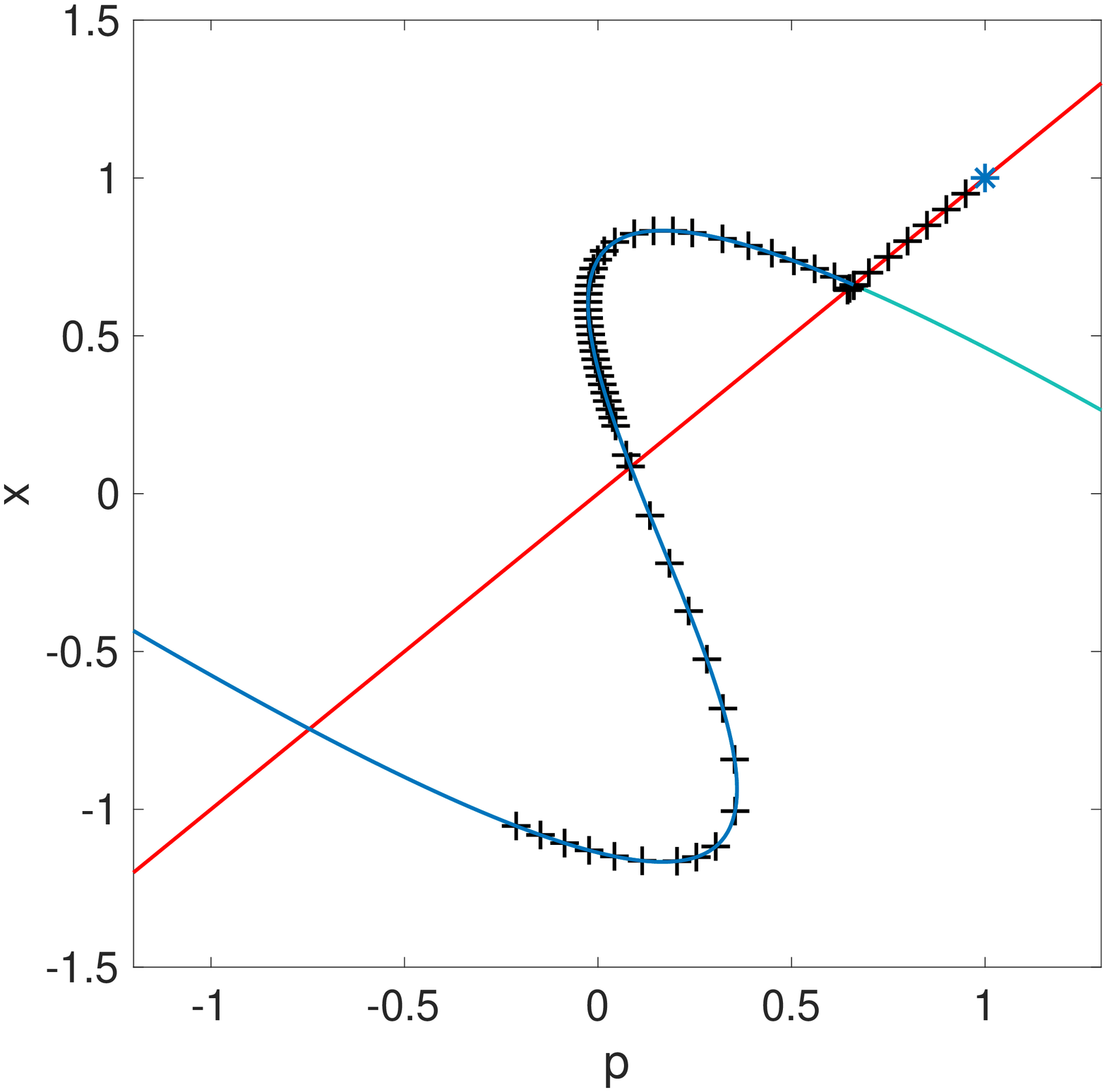}
	\includegraphics[width=0.45\linewidth]{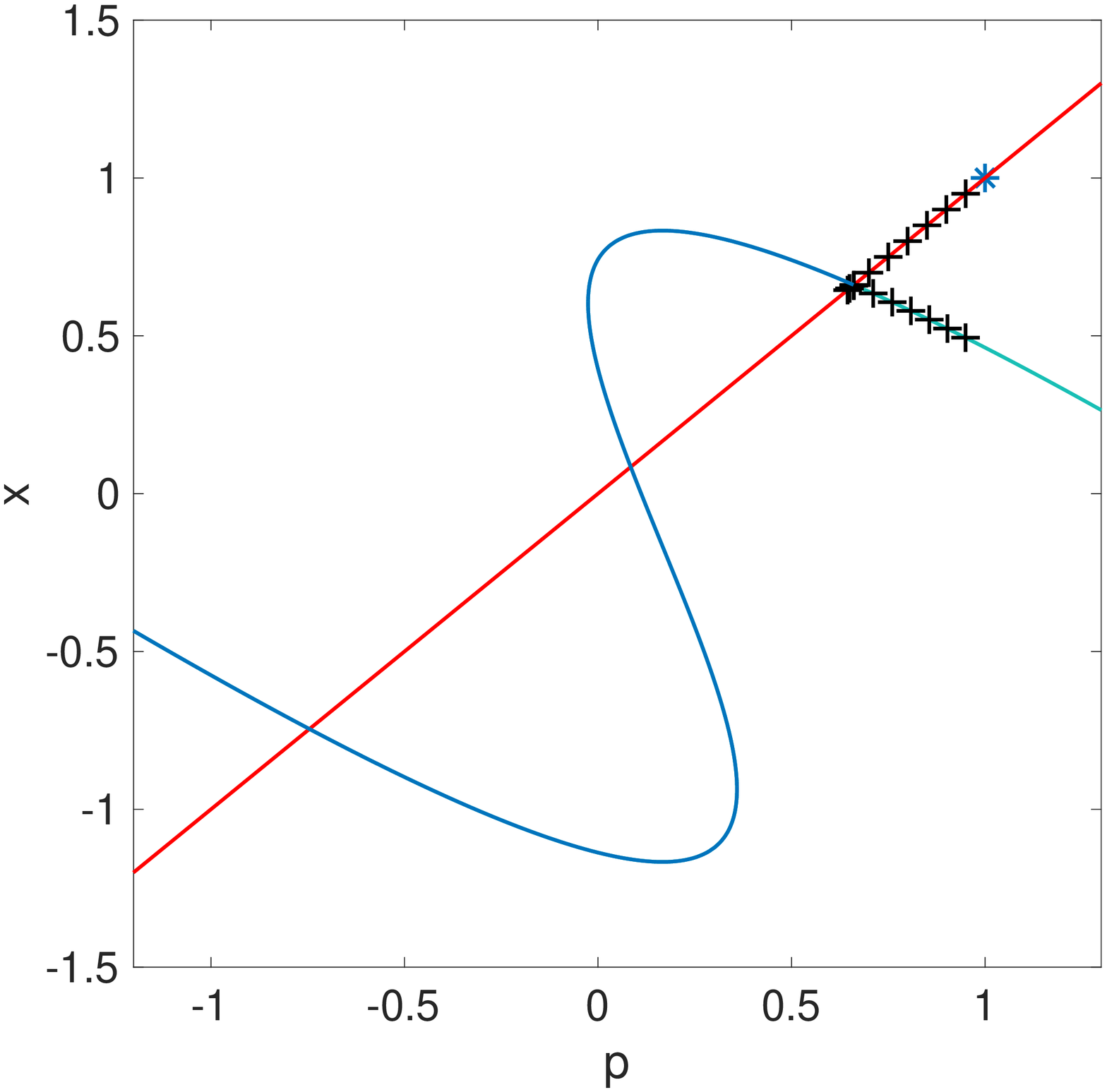}
	\caption{Solution behavior of (\ref{ex23}) with diagonal (red) and non-diagonal (blue) branches.}		
	\label{Fig:ex23}
\end{figure}

\noindent	{\bf Example 3:} We consider the following equation which is used in \cite{wang2016singularity} as an example of universal unfoldings of singularities of topological codimension two:
\begin{equation}
F(x,p)=(x-p)^2+(\frac{1}{3}-2(x+p)+(x+p)^3)(x-p).\label{ex23}
\end{equation}
We used the AHTBD method to track the solution branch starting from $(1,1)$, which is shown as the blue point in Fig. \ref{Fig:ex23} and is tracked along two different solution branches after the first bifurcation point.
We also compared the traditional homotopy tracking with the AHTBD method on two branches: diagonal and non-diagonal (red and blue, respectively, in Fig. \ref{Fig:ex23}). In Table \ref{tab:ex3_compare}, we tracked from $(1,1)$ with $h=-0.05$  until $p<-0.03$. When tracking along the non-diagonal branch, we encountered turning points where the AHTBD method works well. However, for the traditional method, we have to switch the tracking parameter from $p$ to $x$ in order to ensure the tracking process follows the correct direction.
\begin{table}
	\begin{center}
		\begin{tabular}{|c|c|c|c|c|c|}
			\hline
			\multicolumn{2}{|c|}{\multirow{2}{*}{}}
			&\multicolumn{2}{|c|}{Diagonal branch}&\multicolumn{2}{|c|}{Non-diagonal branch}\\\cline{3-6}
			\multicolumn{2}{|c|}{\multirow{2}{*}{}}	& Trial-and-error   & AHTBD & Trial-and-error   & AHTBD\\
			\hline
			\multicolumn{2}{|c|}{\# of steps}  & 139 & 27  &273  & 60  \\ \hline
			{\multirow{2}{*}{1st bifurcation}}& x& $0.6612$  & $0.6609$  &  $0.6612$  & $0.6609$\\\cline{2-6}
			&p	& $0.6612$  & $0.6609$&  $0.6612$  & $0.6609$\\\hline
			{\multirow{2}{*}{2nd bifurcation}}& x& $0.0846$  & $0.0873$  &  $0.0852$  & $0.0860$\\\cline{2-6}
			&p	& $0.0846$  & $0.0837$&  $0.0843$  & $0.0843$\\\hline
		\end{tabular}
	\end{center}
	\caption{Comparisons between AHTBD and traditional homotopy tracking methods for (\ref{ex23}) along two branches.}
	\label{tab:ex3_compare}
\end{table}

\subsection{An example of nonlinear PDEs}
\begin{figure}[th]
	\centering
	{		\includegraphics[width=\linewidth]{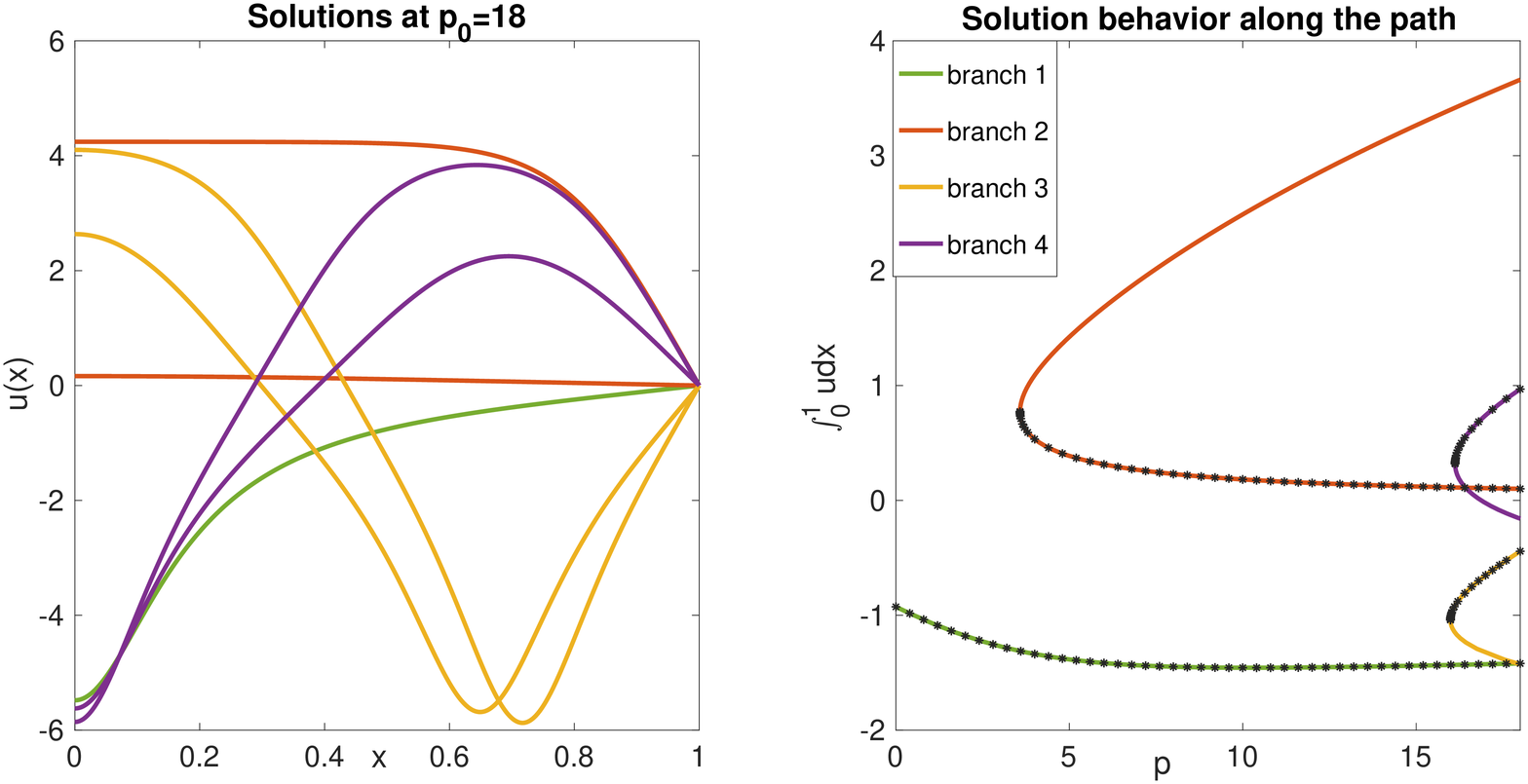}
		\caption{{\bf Left:} Non-trivial solutions of (\ref{1dPDE}) at $p_0=18$; {\bf Right:} Solution behavior of (\ref{1dPDE}) obtained by both the AHTBD method(solid line)  and the traditional method(asterisk).}
		\label{fig:pdeex1}}
\end{figure}

{We compared the AHTBD method with the trial-and-error tracking method on the following nonlinear differential equation:
	\begin{equation*}\left\{
	\begin{aligned}
	&u_{xx} = u^2(u^2-p),\\
	&u_x(0) = 0,	u(1) = 0,
	\end{aligned}\right.\label{1dPDE}
	\end{equation*}
	where $u$ is the solution of differential equation and $p$ is the parameter. There are multiple solutions $u$ for any given parameter $p$, moreover, the number of solutions increases as $p$ goes large. We discretized the differential equation by using the finite difference method and obtained a nonlinear system of polynomial equations. For $p=18$, we solved the discretized nonlinear system by using Newton's method with different initial guesses and obtained seven solutions that is shown in Fig. \ref{fig:pdeex1}. Then we tracked $p$ from $18$ to $0$ with $h=-0.4$ and compared two methods. The stopping criteria for the trial-and-error method is that the stepsize is less than $1e-9$ while it is $p(p-18) >0$ for the AHTBD method. We compared two methods in the tracking steps and running time for the nonlinear system with $360$ grid points in Table \ref{tab:time_step}. The AHTBD method is more efficient to obtain the full solution behaviors for different branches while the traditional trial-and-error tracking method obtains half of branches. 
}
\begin{table}[ht]
	{\footnotesize
		{			\begin{center}
				\begin{tabular}{|c|c|c|c|c|} \hline
					\multirow{2}{*}{Branch No.}&
					\multicolumn{2}{|c|}{ Trial-and-error  }&\multicolumn{2}{|c|}{AHTBD}\\\cline{2-5}
					& Steps& Elapsed time & Steps & Elapsed time\\\hline
					2 & 46 steps&  1.1317s&  31 steps&  1.553s
					\\\hline  3& 57 steps& 2.2933s &  38 steps  &1.7895s
					\\\hline  4 &33 steps & 1.6764s&  15 steps    & 0.7928s
					\\\hline
					5 & 30 steps& 1.6732s &  16 steps  & 0.8745s
					\\\hline
				\end{tabular}
			\end{center}
		}
		\caption{Comparisons between the AHTBD method and the trial-and-error method for the discretized nonlinear system of (\ref{1dPDE}) with $360$ grid points.} \label{tab:time_step}}
\end{table}

\section{Application to a system of nonlinear PDEs}
We apply the AHTBD method to a system of nonlinear PDEs to model two species:
consider a competition between two species that are ecologically identical except in
their dispersal mechanisms. Let $u = u(x)$, $v = v(x)$ denote the densities of two competing species at
location $x$. Then the study of the interaction between a resident phenotype ($u$) with an invader phenotype ($v$) can be modeled by the following system:
\begin{equation}\label{eq:SC}
\left\{ \begin{array}{rcl}
\nabla(d\nabla u-\alpha u \nabla m) &=-u(m-u) & \text{ in }\Omega, \\
\nabla(d\nabla v-\beta v \nabla m) &=-v(m-u ) & \text{ in }\Omega , \\
d \frac{\partial u}{\partial n}  - \alpha u \frac{\partial m}{\partial n} &= d\frac{\partial v}{\partial n } - \beta v \frac{\partial m}{\partial n} = 0 & \text{ on }\partial\Omega. \\
\end{array} \right.
\end{equation}
Here $m(x)$ is the per-capita growth rate, which represents the same resources that two species are competing for. To reflect the heterogeneity of the environment, we assume that $m(x)$
is a nonconstant function to  reflect the quality and quantity of
resources available at the location $x$.
In Eq.~(\ref{eq:SC}), $d$ is two species' common random dispersal rates, and $\alpha, \beta$ are their rates of directed movement upward along the resource gradient. The boundary condition is of a no-flux type, i.e., there is no net movement across the boundary. The solution behavior of this model has been studied well in \cite{CHL,HamLou,HaoLamLou,LamLou}: when $\alpha=\beta$, two species co-exist, $u=v$. Bifurcation, the so-called evolutionarily stable strategy (ESS), happens on the diagonal $\alpha=\beta$, and the behavior of the solution near the bifurcation point is described in \cite{CHL,HamLou,LamLou}.
In reality, it is interesting to find out what happens for the bifurcation branch away
from the bifurcation point, and this is where the numerical
computation is needed: to find the population densities $u$ and $v$ as $\alpha$ and $\beta$ moves far away from the ESS. Given $m(x)=1+x$, a unique positive solution of $u$ is defined by \eqref{eq:SC}, namely, $\tilde u = \tilde{u}(d,\alpha)$.	By standard theory, if some rare population $v$ is introduced into the resident population $u$ at equilibrium
(i.e. $u \equiv \tilde u$), then the initial (exponential) growth rate of the population of $v$ is given by
$\lambda$, where $\lambda  = \lambda(\alpha,\beta; d)$ is the principal eigenvalue of the problem
\begin{equation}\label{eq:1.4}
\left\{
\begin{array}{ll}
\nabla \cdot (d \nabla \varphi - \beta \varphi \nabla m) + (m - \tilde{u}(d,\alpha) ) \varphi = \lambda \varphi&\text{ in }\Omega,\\
d \frac{\partial \varphi}{\partial n} - \beta \varphi \frac{\partial m}{\partial n} = 0 &\text{ on }\partial\Omega,
\end{array}
\right.
\end{equation}
where the positive principal eigenfunction $\varphi  = \varphi(\alpha,\beta; d)$ is uniquely determined by the normalization
\begin{equation}\label{eq:1.5}
\int_\Omega \varphi (\alpha,\beta;d)=1.
\end{equation}
In particular, when $ \alpha=\beta$, we have $\varphi(\alpha,\alpha;d) = \tilde{u}$ and $\lambda(\alpha,\alpha;d) \equiv 0$ for any $d, \alpha$ which implies that two species $u$ and $v$ are identical when $\alpha = \beta$.

\begin{figure}[th]
	\centering
	\includegraphics[width=0.6\linewidth]{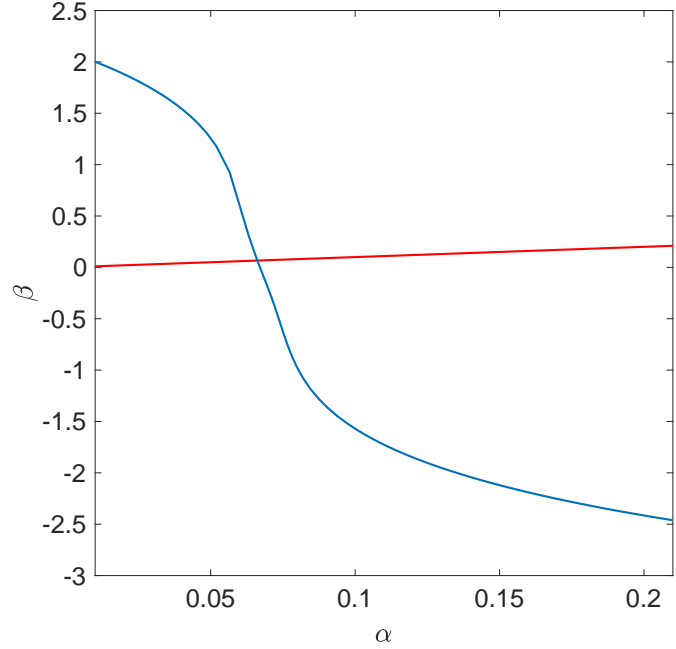}
	\caption{Diagram of $\alpha$-$\beta$ by tracking $\mathbf{F}(\beta,\mathbf{u},\mathbf{v};\alpha)=0$ with respect to $\alpha$.  }\label{Fig:SB}
\end{figure}

When we couple \eqref{eq:SC} and \eqref{eq:1.5} together and discretize the system by using the finite difference method, we have the following coupled system:
\begin{equation}
\label{eq:DC}
\mathbf{F}(\beta,\mathbf{u},\mathbf{v};\alpha):=
\left(\begin{matrix}
\frac{2d}{h^2}u_2-(\frac{2d}{h^2}+\frac{2\alpha}{h}+\frac{\alpha^2}{d})u_1+u_1(m_1-u_1 )\\
\frac{d}{h^2}(u_{i+1}-2u_i+u_{i-1})- \frac{\alpha}{2h}(u_{i+1}-u_{i-1}) +u_i(m_i-u_i )\\
(-\frac{2d}{h^2}+\frac{2\alpha}{h}-\frac{\alpha^2}{d})u_N+ \frac{2d}{h^2}u_{N-1}+u_N(m_N-u_N)\\
\frac{2d}{h^2}v_2-(\frac{2d}{h^2}+\frac{2\beta}{h}+\frac{\beta^2}{d})v_1+v_1(m_1-u_1)\\
\frac{d}{h^2}(v_{i+1}-2v_i+v_{i-1})- \frac{\beta}{2h}(v_{i+1}-v_{i-1}) +v_i(m_i-u_i)\\
(-\frac{2d}{h^2}+\frac{2\beta}{h}-\frac{\beta^2}{d})v_N+ \frac{2d}{h^2}v_{N-1}+v_N(m_N-u_N)\\
(\frac{v_1}{2}+v_2+\cdots+v_{N-1}+\frac{v_N}{2})h-1
\end{matrix}\right)=0.
\end{equation}
For any given $\alpha_0$, $\mathbf{u}_0$ is solved by the discretization of \eqref{eq:SC}. Then $\mathbf{u}_0$, $\beta_0=\alpha_0$, $\mathbf{v}_0=\frac{\mathbf{u}_0}{\int_\Omega \mathbf{u}_0}$ is a solution of $\mathbf{F}(\beta,\mathbf{u},\mathbf{v};\alpha)=0$. Given initial values $(\beta_0,\mathbf{u}_0,\mathbf{v}_0,\alpha_0)$, we track along the diagonal branch $\alpha=\beta$ using $\alpha$ as a
parameter. For our choice of $m(x)$, there is only one bifurcation. We applied the AHTBD method to track $\mathbf{F}(\beta,\mathbf{u},\mathbf{v};\alpha)=0$, which is shown in Fig. \ref{Fig:SB} by starting with $\alpha_0=0.01$ and ending with $\alpha_0>0.3$.
We also compared the AHTBD method with the traditional trial-and-error  tracking method in Tables \ref{tab:different_stepsize} \& \ref{tab:time}  and demonstrated that the AHTBD
method is faster than the traditional homotopy tracking
method for the nonlinear PDE example.  

\begin{table}[ht]
	{		\begin{center}
			\begin{tabular}{|c|c|c|c|c|} \hline
				\multirow{2}{*}{$h$}&
				\multicolumn{2}{|c|}{Diagonal branch}&\multicolumn{2}{|c|}{Non-diagonal branch}\\\cline{2-5}
				& Trial-and-error  & AHTBD &Trial-and-error   & AHTBD\\
				\hline
				0.01 & 88 steps(42.2808s)  & 25 steps(15.5970s)& 88 steps(48.5956s)& 26 steps(16.9497s)\\
				\hline 0.02&  70 steps(33.2463s)  & 16 steps(10.8033s) &  70 steps(40.5963s) & 15 steps(10.2496s)\\\hline
			\end{tabular}
	\end{center}}
	\caption{Comparison between the AHTBD method and the traditional trial-and-error  tracking  with different step-sizes for $h$ (the number of grid points $N=320$).} \label{tab:different_stepsize}
\end{table}

\begin{table}[ht]
	{\footnotesize
		{		\begin{center}
				\begin{tabular}{|c|c|c|c|c|} \hline
					\multirow{2}{*}{$N$}&
					\multicolumn{2}{|c|}{Diagonal branch}&\multicolumn{2}{|c|}{Non-diagonal lower branch}\\\cline{2-5}
					& Trial-and-error tracking & AHTBD& Trial-and-error tracking& AHTBD\\
					\hline 80 & 85 steps(5.7142s)  & 28 steps(2.7299s) &  85 steps(6.1369s)  & 26 steps(2.4906s)\\\hline
					160  & 96 steps(17.9011s) & 29 steps(6.2515s)& 96 steps(19.1336s)  & 28 steps(6.6682s)
					\\\hline
					320  &  88 steps(42.2808s)  & 25 steps(15.5970s)& 88 steps(48.5956s)& 26 steps(16.9497s)\\
					\hline			\end{tabular}
			\end{center}
	}}
	\caption{Comparison between the AHTBD method and the traditional trial-and-error  tracking for number of grid points $N$ (the step-size is $h=0.01$).} \label{tab:time}
\end{table}

\section{Conclusions}
\label{sec:conclusions}
We developed an adaptive homotopy tracking method to compute bifurcations for large-scale nonlinear parametric systems. This new algorithm is designed
for computing bifurcation points and solutions on different branches
through the bifurcations via adaptive tracking, the Puiseux interpolation, and the inflation process. Furthermore, an augmented system is introduced to compute the adaptive parameter step-size while the inflation technique is backed up when the augmented system becomes singular. We also employ the Puiseux series expansion to interpolate bifurcation points, and different bifurcation branches are approximated based on computing the tangent cone structure of the bifurcation point. Several
numerical examples for both polynomial systems and nonlinear systems of PDEs verify the efficiency of this new method through comparison with the traditional homotopy continuation method.
{There are still some numerical challenges for  the adaptive homotopy tracking method developed in this paper. For example, it would become challenging and might fail when we deal with  a cluster of bifurcations. Moreover, the efficient and accurate eigenvalue solver is required in this adaptive tracking process. Thus inexact approximations of eigen data and inaccurate solution points might also affect the numerical performance. We will explore these challenges more carefully in the future.}

\bibliographystyle{unsrt}  


\end{document}